\documentclass [10pt] {article}

 \usepackage [english] {babel}
 \usepackage {amssymb}
 \usepackage {amsmath}
 \usepackage {amsthm}
 \usepackage {xypic}

 \def \^{\hat}
 \def \"{\check}
 \def \${\sim}
 \def \_{\underline}
 \def \={\overline}
 \def \<{\langle}
 \def \>{\rangle}
 \def \[{\mathopen{\pmb{\boldsymbol\langle}}}
 \def \]{\mathclose{\pmb{\boldsymbol\rangle}}}
 \def \({\mathopen{\pmb{\boldsymbol(}}}
 \def \){\mathclose{\pmb{\boldsymbol)}}}
 \def \|{\mathbin|}

 \def \:{\colon}
 \def \?{{-}}

 \def \N{{\boldsymbol{\mathsf N}}}
 \def \R{{\boldsymbol{\mathsf R}}}
 \def \k{{\boldsymbol k}}

 \def \kMod {{\text{\bf$\k$-Mod}}}

 \DeclareMathOperator {\Hom} {Hom}
 \DeclareMathOperator {\Aut} {Aut}
 \DeclareMathOperator {\alt} {alt}
 \DeclareMathOperator {\sgn} {sgn}

 \newcommand* {\head} [1]
 {\subsubsection * {#1}}

 \newcommand* {\subhead} [1]
 {\addvspace\medskipamount \noindent {\bf\itshape #1\/}}

 \newenvironment* {claim} [1] []
 {\begin{trivlist}\item [\hskip\labelsep {\bf #1}] \it}
 {\end{trivlist} }

 \newenvironment* {demo} [1] []
 {\begin{trivlist}\item [\hskip\labelsep {\it #1}] }
 {\end{trivlist} }

\begin {document}

 \title {\Large\bf
         Commutative algebras \\
         and representations
         of the category of finite sets}

 \author {\normalsize\rm
          S.~S.~Podkorytov}

 \date {}

 \maketitle

 \begin {abstract} \noindent
 We prove that
 two finite-dimensional commutative algebras
 over an algebraically closed field
 are isomorphic if and only if
 they give rise to isomorphic representations
 of the category of finite sets and surjective maps.
 \end {abstract}


 Let $\Omega$ be the category
 whose objects are the sets $\_n=\{1,\dotsc,n\}$, $n\in\N$
 ($=\{0,1,\dotsc\}$), and
 whose morphisms are surjective maps.
 Let $\k$ be a field.
 For a [commutative] algebra $A$ without unity [over $\k$],
 let us define a functor $L_A\:\Omega\to\kMod$
 (a ``representation of $\Omega$'').
 For $n\in\N$,
 set $L_A(\_n)=A^{\otimes n}$.
 For a morphism $h\:\_m\to\_n$,
 set $L_A(h)\:x_1\otimes\dotso\otimes x_m\mapsto
 y_1\otimes\dotso\otimes y_n$, where
 $$
 y_j=\prod_{i\in h^{-1}(j)}x_i.
 $$
 The functor $L_A$ is a variant of the Loday functor
 \cite{Wikipedia}.

 \begin {claim} [Theorem.]
 Let the field $\k$ be algebraically closed.
 Let $A$ and $B$ be finite-dimensional algebras without unity.
 Suppose that the functors $L_A$ and $L_B$ are isomorphic.
 Then the algebras $A$ and $B$ are isomorphic.
 \end {claim}

 We do not know whether
 the assertion is true for infinite-dimensional algebras.
 It is false for the field $\R$
 (which is not algebraically closed).
 Indeed,
 take the non-isomorphic algebras
 $A=\R[X]/(X^2-1)$ and
 $B=\R[Y]/(Y^2+1)$.
 We have the bases
 $\{X^e\}_{e=0,1}$ ($=\{1,X\}$) in $A$ and
 $\{Y^e\}_{e=0,1}$ in $B$.
 The linear maps $s_n\:A^{\otimes n}\to B^{\otimes n}$,
 $$
 X^{e_1}\otimes\dotso\otimes X^{e_n}\mapsto
 k_{e_1+\dotso+e_n}Y^{e_1}\otimes\dotso\otimes Y^{e_n}, \qquad
 e_1,\dotsc,e_n=0,1,
 $$
 where $k_m=(-1)^{[m/2]}$,
 form a functor isomorphism $s\:L_A\to L_B$.


 \head {1. Preliminaries}


 \subhead {Algebra of polynomials.}
 If $V$ is a vector space [over $\k$],
 then the symmetric group $\Sigma_n=\Aut\_n$ acts [from the
 left] on $V^{\otimes n}$  by the rule
 $g(v_1\otimes\dotso\otimes v_n)=
 v_{g^{-1}(1)}\otimes\dotso\otimes v_{g^{-1}(n)}$.
 The symmetric powers $S^n(V)=(V^{\otimes n})_{\Sigma_n}$ form
 the symmetric algebra
 $$
 S(V)=\bigoplus_{n=0}^\infty S^n(V)
 $$
 with the multiplication induced by the tensor one:
 $\=x\,\=y=\={x\otimes y}$,
 $x\in V^{\otimes m}$, $y\in V^{\otimes n}$
 (the bar denotes the projection $V^{\otimes n}\to S^n(V)$).

 For a vector space $U$,
 put $\k[U]=S(U^*)$.
 For $u\in U$,
 there is the evaluation map $\k[U]\to\k$, $f\mapsto f(u)$,
 which is the algebra homomorphism defined by the condition
 $v(u)=\<v,u\>$ for $v\in U^*=S^1(U^*)\subseteq\k[U]$.
 For a polynomial $f\in\k[U]$ and a set $X\subseteq U$,
 there is the function $f\|X\:X\to\k$, $u\mapsto f(u)$.
 An ideal $P\subseteq\k[U]$ determines the set
 $$
 Z(P)=\{\,u:\text{$f(u)=0$ for all $f\in P$}\,\}\subseteq U.
 $$


 \subhead {Symmetric tensors, isomorphism $\theta$.}
 Put $D^n(U)=(U^{\otimes n})^{\Sigma_n}$,
 $$
 \^D(U)=\prod_{n=0}^\infty D^n(U).
 $$
 The pairing
 \begin{equation} \label{asym}
 \[\?,\?\]\:
 (U^*)^{\otimes n}\times U^{\otimes n}\to\k,
 \end{equation}
 $\[v_1\otimes\dotso\otimes v_n,u_1\otimes\dotso\otimes u_n\]=
 \<v_1,u_1\>\dotso\<v_n,u_n\>$,
 induces the pairing
 \begin{equation} \label{sym}
 \(\?,\?\)\:
 S^n(U^*)\times D^n(U)\to\k,
 \end{equation}
 $\(\=z,w\)=\[z,w\]$, where
 $w\in D^n(U)\subseteq U^{\otimes n}$,
 $z\in(U^*)^{\otimes n}$.
 Summing over $n\in\N$,
 we get a pairing
 $$
 \(\?,\?\)\:
 \k[U]\times\^D(U)\to\k.
 $$
 We have the linear map
 $$
 \theta\:\^D(U)\to\k[U]^*, \quad
 \<\theta(W),f\>=\(f,W\).
 $$
 If $U$ is finite-dimensional, then
 the pairings \eqref{asym} and \eqref{sym} are perfect
 and $\theta$ is an isomorphism.


 \subhead {Functor $T_A$.}
 Let $\Sigma\subseteq\Omega$ be the subcategory of
 isomorphisms.
 We have $\Sigma=\Sigma_0\sqcup\Sigma_1\sqcup\dotso$.
 For a vector space $A$,
 we have the functor
 $T_A\:\Sigma\to\kMod$, $T_A(\_n)=A^{\otimes n}$
 (with the ordinary action of $\Sigma_n$).
 If $A$ is an algebra without unity,
 then $T_A=L_A|_\Sigma$.


 \subhead {Kronecker product, isomorphism $\kappa$.}
 If a group $G$ acts on vector spaces $X$ and $Y$,
 then it acts on $\Hom(X,Y)$ by the rule
 $(gt)(x)=g(t(g^{-1}x))$.
 We have $\Hom(X,Y)^G=\Hom_G(X,Y)$.

 Let $A$ and $B$ be vector spaces.
 The Kronecker product
 $\Hom(A,B)^{\otimes n}\to\Hom(A^{\otimes n},B^{\otimes n})$,
 $w\mapsto[w]$ (a notation),
 preserves the action of $\Sigma_n$ and
 thus induces a linear map
 $D^n(B^A)\to\Hom_{\Sigma_n}(A^{\otimes n},B^{\otimes n})$
 (from now on, $B^A=\Hom(A,B)$).
 Since
 $$
 \Hom_\Sigma(T_A,T_B)=\prod_{n=0}^\infty
 \Hom_{\Sigma_n}(A^{\otimes n},B^{\otimes n}),
 $$
 these maps form a linear map
 $$
 \kappa\:\^D(B^A)\to\Hom_\Sigma(T_A,T_B).
 $$
 If $A$ and $B$ are finite-dimensional,
 then $\kappa$ is an isomorphism.


 \subhead {Morphisms $T_A\to T_B$ and functionals on $\k[B^A]$,
           isomorphism $\xi$.}
 For finite-dimensional vector spaces $A$ and $B$ we have the
 isomorphism $\xi$
 that fits in the commutative diagram
 $$
 \xymatrix {
 & &
 \Hom_\Sigma(T_A,T_B)
 \ar[dd]^-{\xi} \\
 \^D(B^A)
 \ar[rru]^-{\kappa}
 \ar[rrd]_-{\theta} & & \\
 & &
 \k[B^A]^*.
 }
 $$

 \noindent {\it Example.\/}
 A linear map $u\:A\to B$ induces the functor morphism
 $T_u\:T_A\to T_B$, $(T_u)_n=u^{\otimes n}$.
 Then $\<\xi(T_u),f\>=f(u)$, $f\in\k[B^A]$.


 \subhead {Antisymmetrization.}
 For a vector space $V$,
 we have the operator $\alt_n\:V^{\otimes n}\to V^{\otimes n}$,
 $$
 \alt_n(w)=\sum_{g\in\Sigma_n}\sgn g\ gw.
 $$


 \head {2. The determinant}


 Let $A$ and $B$ be vector spaces of equal finite dimension
 $m$.
 Put $U=B^A$.
 Choose bases $e_1,\dotsc,e_m\in A$ and $f_1,\dotsc,f_m\in B$.
 Put
 $$
 E=\alt_m(e_1\otimes\dotso\otimes e_m)\in A^{\otimes m}, \quad
 F=\alt_m(f_1\otimes\dotso\otimes f_m)\in B^{\otimes m}.
 $$
 We have the bases
 $\"f^1,\dotsc,\"f^m\in B^*$, $\<\"f^j,f_i\>=\delta_i^j$
 ($\delta_i^j$ is the Kronecker delta) and
 $l_i^j\in U^*$, $i,j=1,\dotsc,m$,
 $\<l_i^j,u\>=\<\"f^j,u(e_i)\>$.
 Put
 $$
 H=\sum_{g\in\Sigma_m}\sgn g\ %
 l_{g^{-1}(1)}^1\otimes\dotsc\otimes l_{g^{-1}(m)}^m
 \in(U^*)^{\otimes m}.
 $$
 Then $\=H\in\k[U]$ is the determinant,
 so
 \begin{equation} \label{H=det}
 \=H(u)=\det u, \qquad
 u\in U.
 \end{equation}

 We have $\[\"f^1\otimes\dotso\otimes\"f^m,[v](E)\]=\[H,v\]$,
 $v\in U^{\otimes m}$.
 Hence
 \begin{equation} \label{wE-asym}
 \[(\"f^1\otimes\dotso\otimes\"f^m)^{\otimes r},
   [w](E^{\otimes r})\]=\[H^{\otimes r},w\], \qquad
 w\in U^{\otimes mr}, \quad
 r\in\N.
 \end{equation}
 For $w\in D^{mr}(U)$,
 we have
 \begin{equation} \label{wE-sym}
 [w](E^{\otimes r})=\(\=H^r,w\)F^{\otimes r}.
 \end{equation}
 Indeed,
 $E^{\otimes r}$ belongs to the image of
 $\alt_m^{\otimes r}\:A^{\otimes mr}\to A^{\otimes mr}$.
 The image of
 $\alt_m^{\otimes r}\:B^{\otimes mr}\to B^{\otimes mr}$ is
 generated by $F^{\otimes r}$ since
 the image of $\alt_m\:B^{\otimes m}\to B^{\otimes m}$ is
 generated by $F$.
 The map $[w]\:A^{\otimes mr}\to B^{\otimes mr}$ preserves the
 action of $\Sigma_{mr}$ and
 thus commutes with $\alt_m^{\otimes r}$.
 Therefore,
 $[w](E^{\otimes r})=tF^{\otimes r}$ for some $t\in\k$.
 From \eqref{wE-asym},
 we get $t=\[H^{\otimes r},w\]=\(\=H^r,w\)$.

 For a morphism $s\:T_A\to T_B$,
 we have
 \begin{equation} \label{sE}
 s_{mr}(E^{\otimes r})=\<\xi(s),\=H^r\>F^{\otimes r}, \qquad
 r\in\N.
 \end{equation}
 This follows from \eqref{wE-sym}:
 if $s=\kappa(W)$, $W\in\^D(U)$, then
 $s_{mr}=[W_{mr}]$ and
 $\<\xi(s),\=H^r\>=\(\=H^r,W_{mr}\)$.


 \head {3. Homomorphisms $A\to B$ and morphisms $L_A\to L_B$}


 Let $A$ and $B$ be finite-dimensional algebras without unity.
 Put $U=B^A$.


 \subhead {Multiplicativity ideal.}
 Take $x,y\in A$ and $p\in B^*$.
 We have the linear form $I_{x,y}^p\in U^*$,
 $$
 \<I_{x,y}^p,u\>=\<p,u(xy)\>, \qquad
 u\in U
 $$
 (the multiplication in $A$ is used) and
 the tensor $J_{x,y}^p\in(U^*)^{\otimes2}$,
 $$
 \[J_{x,y}^p,u\otimes v\]=\<p,u(x)v(y)\>, \qquad
 u,v\in U
 $$
 (the multiplication in $B$ is used).
 Put
 $$
 g_{x,y}^p=\={J_{x,y}^p}-I_{x,y}^p\in\k[U].
 $$
 We have
 $$
 g_{x,y}^p(u)=\<p,u(x)u(y)-u(xy)\>, \qquad
 u\in U.
 $$
 Let $M\subseteq\k[U]$ be the ideal generated by the
 polynomials $g_{x,y}^p$, $x,y\in A$, $p\in B^*$.

 \begin {claim} [Lemma 1.]
 The set $Z(M)\subseteq U$ coincides with
 the set of algebra homomorphisms $A\to B$.
 \qed
 \end {claim}

 Note that $\Hom_\Omega(L_A,L_B)\subseteq\Hom_\Sigma(T_A,T_B)$.

 \begin {claim} [Lemma 2.]
 Let $s\in\Hom_\Sigma(T_A,T_B)$.
 Then the conditions
 $s\in\Hom_\Omega(L_A,L_B)$ and
 $\xi(s)\perp M$
 are equivalent.
 \end {claim}

 Thus we establish an isomorphism
 $\Hom_\Omega(L_A,L_B)\to(\k[U]/M)^*$.

 \begin {demo} [Proof.]
 For $n\in\N$,
 define the morphism $\tau_n\:\_{n+2}\to\_{n+1}$ by the rules
 $1\mapsto1$ and
 $i\mapsto i-1$, $i>1$.
 The category $\Omega$ is obtained from $\Sigma$
 by adjunction of the morphisms $\tau_n$.
 Therefore,
 the condition $s\in\Hom_\Omega(L_A,L_B)$ is equivalent to
 commutativity of the diagrams
 $$
 \xymatrix {
 A^{\otimes(n+2)}
 \ar[rr]^-{s_{n+2}}
 \ar[d]_-{L_A(\tau_n)} & &
 B^{\otimes(n+2)}
 \ar[d]^-{L_B(\tau_n)} \\
 A^{\otimes(n+1)}
 \ar[rr]^-{s_{n+1}} & &
 B^{\otimes(n+1)},
 }
 $$
 $n\in\N$.
 Consider the discrepancy
 $$
 r_n=L_B(\tau_n)\circ s_{n+2}-s_{n+1}\circ L_A(\tau_n)\:
 A^{\otimes(n+2)}\to B^{\otimes(n+1)}.
 $$

 For $z\in A$, $q\in B^*$,
 we have the linear form $l_z^q\in U^*$,
 $\<l_z^q,u\>=\<q,u(z)\>$, $u\in U$.
 These forms generate $U^*$.
 For $n\in\N$, $x,y,z_1,\dotsc,z_n\in A$,
 $p,q_1,\dotsc,q_n\in B^*$,
 put
 $$
 G_{x,y,z_1,\dotsc,z_n}^{\ p,\ q_1,\dotsc,q_n}=
 g_{x,y}^pl_{z_1}^{q_1}\dotso l_{z_n}^{q_n}\in\k[U].
 $$
 These polynomials linearly generate $M$.
 Therefore, it suffices to prove that
 \begin{equation*}
 \[p^\$,r_n(x^\$)\]=
 \<\xi(s),G_{x,y,z_1,\dotsc,z_n}^{\ p,\ q_1,\dotsc,q_n}\>,
 \end{equation*}
 where
 $x^\$=x\otimes y\otimes z_1\otimes\dotso\otimes z_n$,
 $p^\$=p\otimes q_1\otimes\dotso\otimes q_n$.

 We have
 \begin{gather*}
 \[p^\$,[w_1](L_A(\tau_n)(x^\$))\]=
 \[I_{x,y}^p\otimes l^\$,w_1\], \qquad
 w_1\in U^{\otimes(n+1)}, \\
 \[p^\$,L_B(\tau_n)([w_2](x^\$))\]=
 \[J_{x,y}^p\otimes l^\$,w_2\], \qquad
 w_2\in U^{\otimes(n+2)},
 \end{gather*}
 where $l^\$=l_{z_1}^{q_1}\otimes\dotso\otimes l_{z_n}^{q_n}$
 (direct check).
 By construction,
 $$
 G_{x,y,z_1,\dotsc,z_n}^{\ p,\ q_1,\dotsc,q_n}=
 \={J_{x,y}^p\otimes l^\$}-
 \={I_{x,y}^p\otimes l^\$}.
 $$

 We have $s=\kappa(W)$ for some sequence $W\in\^D(U)$,
 so $s_n=[W_n]$.
 We have
 \begin{multline*}
 \[p^\$,r_n(x^\$)\]=
 \[p^\$,L_B(\tau_n)([W_{n+2}](x^\$))\]-
 \[p^\$,[W_{n+1}](L_A(\tau_n)(x^\$))\]= \\ =
 \[J_{x,y}^p\otimes l^\$,W_{n+2}\]-
 \[I_{x,y}^p\otimes l^\$,W_{n+1}\]= \\ =
 \<\theta(W),G_{x,y,z_1,\dotsc,z_n}^{\ p,\ q_1,\dotsc,q_n}\>=
 \<\xi(s),G_{x,y,z_1,\dotsc,z_n}^{\ p,\ q_1,\dotsc,q_n}\>.
 \qed
 \end{multline*}
 \end {demo}


 \subhead {Proof of Theorem.}
 Let $s\:L_A\to L_B$ be a functor isomorphism.
 Then $s_1\:A\to B$ is an isomorphism of vector spaces.
 Put $m=\dim A=\dim B$.
 Choose bases in $A$ and $B$.
 Let the tensors $E$, $F$ and $H$ be as in \S~2.
 We seek an algebra homomorphism $u\:A\to B$ with $\det u\ne0$.
 Assume that there exists no such a homomorphism.
 Then,
 by \eqref{H=det} and Lemma 1,
 $\=H\|Z(M)=0$.
 By Hilbert's Nullstellensatz,
 $\=H^r\in M$ for some $r\in\N$.
 By \eqref{sE} and Lemma 2,
 $s_{mr}(E^{\otimes r})=\<\xi(s),\=H^r\>F^{\otimes r}=0$.
 This is absurd since
 $E^{\otimes r}\ne0$ and
 $s_{mr}$ is an isomorphism of vector spaces.
 \qed



 {\noindent \tt ssp@pdmi.ras.ru}

 {\noindent \tt http://www.pdmi.ras.ru/\~{}ssp}

 \end {document}